\documentclass[12pt]{amsart}

\usepackage{amssymb}

\usepackage{enumerate}

\usepackage{graphicx}

\makeatletter
\@namedef{subjclassname@2010}{%
  \textup{2010} Mathematics Subject Classification}
\makeatother

\newtheorem{theorem}{ Main Theorem}[section]
\newtheorem{thm}{Theorem}[section]

\newtheorem{lem}[thm]{Lemma}

\theoremstyle{definition}

\newtheorem{rem}[thm]{Remark}

\numberwithin{equation}{section}

\begin{document}


\baselineskip=17pt


\title[A sum of cubes equals a sum of 5th powers]{On the Diophantine equation in the form that a sum of cubes equals a sum of quintics }

\author[F. Izadi]{Farzali Izadi}
\address{Farzali Izadi \\
Department of Mathematics \\ Faculty of Science \\ Urmia University \\ Urmia 165-57153, Iran}
\email{f.izadi@urmia.ac.ir}

\author[M. Baghalaghdam]{Mehdi Baghalaghdam}
\address{Mehdi Baghalaghdam \\
Department of Mathematics\\ Faculty of Science \\ Azarbaijan Shahid Madani University\\Tabriz 53751-71379, Iran}
\email{mehdi.baghalaghdam@azaruniv.edu}

\date{}

\begin{abstract}
In this paper, the elliptic curves theory is used for solving the Diophantine equations
$a({X'_1}^{5}+{X'_2}^{5})+ \sum_{i=0}^n a_i {X_{i}} ^{5}=b({Y'_1}^{3}+{Y'_2}^{3})+\sum_{i=0}^m b_i{Y_i}^{3}$,
where $n,m\in \mathbb{N}\cup\{0\}$, and, $a,b\neq0$, $a_i$, $b_i$, are fixed arbitrary rational numbers. We solve  the Diophantine equation for some values of $n$, $m$, $a$, $b$, $a_i$, $b_i$,  and obtain nontrivial integer solutions for each case. By our method, we may find infinitely many nontrivial integer solutions for the Diophantine equation for every  $n$, $m$, $a$, $b$, $a_i$, $b_i$, and show among the other things that how  sums of some $5th$ powers can be written as sums of some cubics.
\end{abstract}



\maketitle

\section{Introduction}
\noindent 
Euler conjectured in $1969$ that the Diophantine equation $A^4+B^4+C^4=D^4$, or more generally $A_1^N+A_2^N+ \cdots +A_{N-1}^N=A_N^N$, ($N\geq4$), has no solution in positive integers ( see \cite{1}). Nearly two centuries later, a computer search (see \cite{6}) found the first counterexample to the general conjecture (for $N=5$):
\\

$27^5+84^5+110^5+133^5=144^5$.
\\

\noindent In $1986$, Noam Elkies, by elliptic curves, found counterexamples for the $N=4$ case (see \cite{2}). His smallest counterexample was:
\\

$2682440^4+15365639^4+18796760^4=20615673^4$.
\\

\noindent The authors in three different papers, used  elliptic curves to solve three Diophantine equations

\begin{equation}\label{1}
{ \sum_{i=1}^n a_ix_{i} ^4= \sum_{j=1}^na_j y_{j}^4 },
\end{equation}
\noindent where $a_i$, and $n\geq3$, are fixed arbitrary integers,\\

\begin{equation}\label{2}
X^4+Y^4= 2U^4+ \sum_{i=1}^nT_i U_{i} ^{\alpha_{i}},
\end{equation}
\noindent where $n, \alpha_i\in \mathbb{N}$, and $T_i$, are appropriate fixed arbitrary rational numbers, and,\\
\begin{equation}\label{3}
\sum_{i=1}^n a_ix_{i} ^6+\sum_{i=1}^m b_iy_{i} ^3= \sum_{i=1}^na_iX_{i}^6\pm\sum_{i=1}^m b_iY_{i} ^3,
\end{equation}
\noindent where $n$, $m$ $\geq 1$ and $a_i$, $b_i$, are fixed arbitrary nonzero integers. (see \cite{3}, \cite{4}, \cite{5})
 \\
 
\noindent In this paper, we are interested in the study of the Diophantine equation:
\\

\begin{equation}\label{4}
a({X'_1}^{5}+{X'_2}^{5})+ \sum_{i=0}^n a_i {X_{i}} ^{5}=b({Y'_1}^{3}+{Y'_2}^{3})+\sum_{i=0}^m b_i{Y_i}^{3},
\end{equation}
\\

\noindent where $n,m\in \mathbb{N}\cup\{0\}$, and, $a,b\neq0$, $a_i$, $b_i$, are fixed arbitrary rational numbers.
\\

\section{our main theorem}

Our main result is the following:
\begin{theorem} Consider the Diophantine equation \eqref{4}.

Let $Y^2=X^3+FX^2+GX+H$, be an elliptic curve
in which the coefficients $F$, $G$, and $H$, are all functions of $a$, $b$, $a_i$, $b_i$ and the other rational parameters $\alpha_i$, $\beta_i$, $t$, $x_1$, $v$, yet to be found later. If the elliptic curve has positive rank, depending on the values of $\alpha_i$, $\beta_i$, $x_1$, (This is done by choosing appropriate arbitrary values for $\alpha_i$($0\leq i \leq n$),
 $\beta_i$($0\leq i \leq m$), $x_1$), the Diophantine equation has infinitely many  integer solutions.
\\

\begin{proof} We solve the Diophantine equation \eqref{4},
if we find rational solutions for the  above Diophantine equation, then by canceling the denominators of $X'_1$, $X'_2$, $X_i$, $Y'_1$, $Y'_2$, $Y_i$, and  by multiplying the both sides of Diophantine equation by  the appropriate value of $M$, we may obtain integer solutions for the Diophantine equation.\\
 Note that if
\noindent $(X'_1,X'_2,X_0, \cdots ,X_n,Y'_1,Y'_2,Y_0, \cdots , Y_m)$,

\noindent is a rational solution for the Diophantine equation  \eqref{1}, then for every  arbitrary rational number $\mu,$

\noindent $(\mu^3X'_1,\mu^3X'_2,\mu^3X_0, \cdots ,\mu^3X_n,\mu^5Y'_1,\mu^5Y'_2,\mu^5Y_0, \cdots , \mu^5Y_m)$

\noindent is a solution for \eqref{1}, too.
\\

\noindent Let: $X'_1=t+x_1$, $X'_2=t-x_1$, $Y'_1=t+v$, $Y'_2=t-v$, $X_i=\alpha_it$, $Y_i=\beta_it$,
where all variables are rational numbers. By substituting these variables in the above Diophantine equation, we get:
 \\

\begin{equation}\label{5}
a(2t^5+10x_1^4t+20x_1^2t^3)+\sum_{i=0}^na_i \alpha_{i} ^5t^5=b(2t^3+6tv^2)+\sum_{i=0}^mb_i \beta_{i} ^3t^3.
\end{equation}
\\

\noindent Then after some simplifications and clearing  the case of $t=0$, we obtain:
\\

\begin{equation}\label{6}
v^2=(\frac{2a+\sum_{i=0}^na_i \alpha_{i} ^5}{6b})t^4+(\frac{20ax_1^2-2b-\sum_{i=0}^mb_i \beta_{i} ^3}{6b})t^2+(\frac{5a}{3b})x_1^4.
\end{equation}
\\

\noindent Now by choosing appropriate arbitrary values for $\alpha_i$ ($0\leq i \leq n$),\\
 $\beta_i$ ($0\leq i \leq m$), $x_1$, such that the rank of the  quartic elliptic curve
\eqref{4} to be positive, and by calculating
$X'_1$, $X'_2$, $X_i$, $Y'_1$, $Y'_2$, $Y_i$, from the relations
\\

$X'_1=t+x_1$, $X'_2=t-x_1$, $X_i=\alpha_it$, $Y'_1=t+v$, $Y'_2=t-v$, $Y_i=\beta_it$,
\\

\noindent some simplifications  and canceling  the  denominators of
$X'_1$, $X'_2$, $X_i$, $Y'_1$, $Y'_2$, $Y_i$,
we obtain infinitely many integer solutions for the  Diophantine equation. The proof of the theorem is complete.
\end{proof}
\end{theorem}

\begin{rem} If in the quartic elliptic curve \eqref{6},
\begin{equation}
Q:=(\frac{5a}{3b})x_1^4,
 \end{equation}
 to be square (It is right for appropriate values of $a$, $b$.), say $q^2$, we  may use the following lemma for transforming this quartic to a cubic  elliptic curve of the form
 \\

$y^2+a_1xy+a_3y=x^3+a_2x^2+a_4x+a_6$, where $a_i\in \mathbb{Q}$. 
\\
\end{rem}
\noindent Then we solve the cubic elliptic curve just obtained of the rank$\geq1$, and get infinitely many solutions for the Diophantine equation \eqref{4}.
\\

\noindent Generally, it is not necessary for $Q$ to be square as we may transform the quartic \eqref{6} to a new quartic in which the constant number is square if the rank of the quartic \eqref{6} is positive. The only important thing is that the rank of the quartic elliptic curve \eqref{6}, be positive for getting infinitely many solutions for \eqref{4}. See the example 1.
\\

\begin{lem} Let K be a field of characteristic not equal to $2$. Consider the equation

$v^2=au^4+bu^3+cu^2+du+q^2$, with $a$, $b$, $c$, $d$ $\in K$.
\\

\noindent Let $x=\frac{2q(v+q)+du}{u^2}$, $y=\frac{4q^2(v+q)+2q(du+cu^2)-(\frac{d^2u^2}{2q})}{u^3}$.
\\

\noindent Define
$a_1=\frac{d}{q}$, $a_2=c-(\frac{d^2}{4q^2})$, $a_3=2qb$, $a_4=-4q^2a$, $a_6=a_2a_4$.
\\

\noindent Then $y^2+a_1xy+a_3y=x^3+a_2x^2+a_4x+a_6$.
\\

\noindent The inverse transformation is
$u=\frac{2q(x+c)-(\frac{d^2}{2q})}{y}$, $v=-q+\frac{u(ux-d)}{2q}$.
\\

\noindent The point $(u, v)=(0,q)$ corresponds to the point $(x, y)=\infty$ and

\noindent $(u, v)=(0,−q)$ corresponds to $(x, y)=(−a_2,a_1a_2-a_3)$. (see \cite{8})
\end{lem}
\begin{rem}
 If in the Diophantine equation \eqref{4}, $n$, $m$ to be odd,\\
  $a_{2k}=a_{2k+1}$, ($0\leq k\leq \frac{n-1}{2}$), and $b_{2k}=b_{2k+1}$, ($0\leq k\leq \frac{m-1}{2}$), we may by second beautiful method transform the Diophantine equation \eqref{4} to another quaratic elliptic curve and then solve it for getting infinitely many nontrivial solutions for the Diophantine equation.
\\

\noindent Now we have 
\begin{equation*}
a({X'_1}^5+{X'_2}^5)+a_0(X_0^5+X_1^5)+a_2(X_2^5+X_3^5)+ \cdots +a_{n-1}(X_{n-1}^5+X_n^5)=
\end{equation*}
\begin{equation*}
b({Y'_1}^3+{Y'_2}^3)+b_0(Y_0^3+Y_1^3)+b_2(Y_2^3+Y_3^3)+ \cdots +b_{m-1}(Y_{m-1}^3+Y_m^3).
\end{equation*}
\\

\noindent Then (after renaming the coffecients, the number of terms, and the variables) the above Diophantine equation (or \eqref{4}) is in the form

\begin{equation}\label{7}
\sum_{i=0}^NA_i( {Z_{i}} ^5+{Z'_{i+1}}^5)=\sum_{i=0}^MB_i( {W_{i}} ^3+{W'_{i+1}}^3).
\end{equation}
\\

\noindent Let $Z_i=t+x_i$, $Z'_{i+1}=t-x_i$, $W_i=t+y_i$, $W'_{i+1}=t-y_i$ ($i\geq0$).
 \\

\noindent By substituting these variables in the above Diophantine equation, we get:
 \\

\begin{equation}\label{8}
\sum_{i=0}^NA_i(2t^5+10x_i^4t+20x_i^2t^3)=\sum_{i=0}^MB_i(2t^3+6ty_i^2).
\end{equation}
\\

\noindent Then after some simplifications and clearing  the case of $t=0$, and letting $y_0=v$, we obtain:
\\

\begin{equation}\label{9}
v^2=(\frac{\sum_{i=0}^N A_i}{3B_0})t^4+(\frac{10\sum_{i=0}^N A_ix_i^2-\sum_{i=0}^M B_i}{3B_0})t^2+(\frac{5\sum_{i=0}^N A_ix_i^4-3\sum_{i=1}^M B_iy_i^2}{3B_0}).
\end{equation}
\\

\noindent Now by choosing appropriate values for $x_i$ ($0 \leq i \leq N$), $y_i$ ($1 \leq i \leq M$), such that the rank of the quartic elliptic curve \eqref{9} to be positive, we obtain infinitely many solutions for the Diophantine equation \eqref{4}. The proof of the theorem is complete.
\end{rem}

\section{Application to Examples}
Now we are going to solve some couple of examples:
\\

\subsection{Example: $X_1^5+X_2^5+X_3^5=Y_1^3+Y_2^3+Y_3^3$}\noindent
\\

\noindent i.e., the sum of $3$ fifth powers can be written as the sum of $3$ cubics.
\\

\noindent Let: $X_1=t+x_1$, $X_2=t-x_1$, $X_3=\alpha t$, $Y_1=t+v$, $Y_2=t-v$, $Y_3=\beta t$.
\\

\noindent Then we get:

\begin{equation}\label{25}
v^2=\frac{2+\alpha^5}{6}t^4+\frac{20x_1^2-2-\beta^3}{6}t^2+(\frac{5}{3})x_1^4.
\end{equation}
 \\

\noindent Note that $(\frac{5}{3})x_1^4$, is not a square. Then we may not use from the above  theorem for transforming the  above quartic to a cubic elliptic curve, but we do this work by another beautiful method. Let us take $x_1=1$, $\alpha=\beta=2$. Then the quartic \eqref{25} becomes
 \\

\begin{equation}\label{26}
 v^2=\frac{17}{3}t^4+\frac{5}{3}t^2+\frac{5}{3}.
\end{equation}
\\

\noindent By searching, we see that the above quartic has two rational points $P_1=(1,3)$, and $P_2=(7,117)$, among others. Let us put $T=t-1$. Then we get
\\

\begin{equation}\label{27}
 v^2=\frac{17}{3}T^4+\frac{68}{3}T^3+\frac{107}{3}T^2+26T+9.
\end{equation}
\\

\noindent Now with the inverse transformation
\begin{equation}\label{28}
 T=\frac{6(X+\frac{107}{3})-\frac{26^2}{6}}{Y},
 \end{equation}

and

\begin{equation}\label{29}
v=-3+\frac{T(XT-26)}{6},
\end{equation}

\noindent the quaratic \eqref{27}, maps to the cubic  elliptic curve

 \begin{equation}\label{30}
 Y^2+\frac{26}{3}XY+136Y=X^3+\frac{152}{9}X^2-204X-\frac{10336}{3}.
 \end{equation}
\\

\noindent The rank of this elliptic curve is $2$ and its generator are the points

\noindent $P_1=(X',Y')=(\frac{-152}{9},\frac{280}{27})$, and $P_2=(X'',Y'')=(\frac{-44}{3},\frac{20}{9})$.
 \\

\noindent To square the left hand of \eqref{30}, let us put $M=Y+\frac{13}{3}X+68$. Then the cubic elliptic curve \eqref{30} transforms to the Weierstrass form
 \\

\begin{equation}\label{31}
 M^2=X^3+\frac{107}{3}X^2+\frac{1156}{3}X+\frac{3536}{3}.
\end{equation}
 \\

\noindent The generators for this new cubic elliptic curve are the two points $G_1=(X',M')=(\frac{-44}{3}, \frac{20}{3})$, and $G_2=(X'',M'')=(\frac{-152}{9},\frac{140}{27})$. Thus we conclude that we could transform the main quartic \eqref{26} to the cubic elliptic curve \eqref{31} of rank equal to $2$.
  \\

\noindent Because of this, the  above cubic elliptic curve has infinitely many rational points and we may obtain infinitely many solutions for the Diophantine equation too. \\
 Since $G_1=(X',M')=(\frac{-44}{3},\frac{20}{3})$, we get $(t,v)=(7,-117)$, that is on the \eqref{26}, by calculating $X_i$, $Y_i$,  from the above relations and after  some simplifications  and canceling  the  denominators of $X_i$, $Y_i$, we obtain a solution for the Diophantine equation

 $X_1^5+X_2^5+X_3^5=Y_1^3+Y_2^3+Y_3^3$ as
 \\

$8^5+6^5+14^5=(-110)^3+124^3+14^3$.
\\

\noindent It is iteresting that we see, $8+6+14=(-110)+124+14$, too.
\\

\noindent Also we have $2G_2=(\frac{373}{36},\frac{-21721}{216})$.
 \\

\noindent By using this new point $2G_2=(\frac{373}{36},\frac{-21721}{216})$, we get $(t,v)=(\frac{11}{47},\frac{2943}{2209})$, on the \eqref{26}, by calculating $X_i$, $Y_i$,  from the above relations and after  some simplifications  and canceling  the  denominators of $X_i$, $Y_i$, we obtain another solution for the Diophantine equation as
 \\

 $128122^5+(-79524)^5+48598^5=359227580^3+(-251874598)^3+107352982^3$.
\\

\noindent By choosing the other points on the elliptic curve such as

\noindent $nG_1$, $nG_2$ ($n=3, 4$, $\cdots$ ) we obtain infinitely many solutions for the  Diophantine equation.

\subsection{ Example: $5.({X'_1}^5+{X'_2}^5)=3.({Y'_1}^3+{Y'_2}^3)$}\noindent
\\

\noindent Let: $X'_1=t+x_1$, $X'_2=t-x_1$, $Y'_1=t+v$, $Y'_2=t-v$.
\\

\noindent Then we get:

\begin{equation}\label{10}
v^2=\frac{5}{9}t^4+\frac{50x_1^2-3}{9}t^2+(\frac{25}{9})x_1^4.
\end{equation}
 \\

\noindent Note that for every arbitrary rational number $x_1$, $(\frac{25}{9})x_1^4$, is a square, then by using the above theorem, we may transform the above quartic to a cubic elliptic curve. Now, we must choose appropriate value for $x_1$, such that the rank of the above quartic to bo positive. Let
 $x_1=1$.
 \\

\noindent Then the quartic \eqref{10} becomes
\\

\begin{equation}\label{862}
v^2=\frac{5}{9}t^4+\frac{47}{9}t^2+(\frac{25}{9})x_1^4.
\end{equation}
\\

\noindent With the inverse transformation
\begin{equation}\label{11}
 t=\frac{\frac{10}{3}(X+\frac{47}{9})}{Y},
 \end{equation}

and

\begin{equation}\label{12}
v=\frac{-5}{3}+\frac{t^2X}{\frac{10}{3}},
\end{equation}

\noindent the quaratic \eqref{862}, maps to the cubic  elliptic curve

 \begin{equation}\label{13}
 Y^2=X^3+\frac{47}{9}X^2-\frac{500}{81}X-\frac{23500}{729}.
 \end{equation}
\\

\noindent The rank of this elliptic curve is $1$ and its generator is the point

\noindent $P=(X,Y)=(\frac{-609566}{164025},\frac{-225298052}{66430125})$. Because of this, the  above elliptic curve has infinitely many rational points and we may obtain infinitely many solutions for the Diophantine equation too.\\
Since $P=(X,Y)=(\frac{-609566}{164025},\frac{-225298052}{66430125})$,
 we get $(t,v)=(\frac{-335475}{226658},\frac{-633289965055}{154121546892})$, by calculating $X'_1$, $X'_2$, $Y'_1$, $Y'_2$, from the above relations and after  some simplifications  and canceling  the  denominators of  $X'_1$, $X'_2$, $Y'_1$, $Y'_2$, we obtain a solution for the Diophantine equation

 $5.({X'_1}^5+{X'_2}^5)=3.({Y'_1}^3+{Y'_2}^3)$ as
 \\

$X'_1=-150939399313320876$,\\

$X'_2=-779731267671365724$,\\

$Y'_1=-812465974773035511128800250760$,\\

$Y'_2=382156766244967634925328109160$.\\
\\

\noindent By choosing the other points on the elliptic curve such as

\noindent $nP$, ($n=2, 3$, $\cdots$ ) we obtain infinitely many solutions for the  Diophantine equation.

\subsection{Example: $n.(X_1^5+X_2^5+X_3^5+X_4^5)=m.(Y_1^3+Y_2^3+Y_3^3+Y_4^3)$ }\noindent
\\

\noindent By letting, $X_1=t+x_1$, $X_2=t-x_1$, $X_3=t+x_2$, $X_4=t-x_2$, $Y_1=t+v$,

\noindent $Y_2=t-v$, $Y_3=t+y_2$, $Y_4=t-y_2$,
\\

\noindent we get:

\begin{equation}\label{14}
v^2=(\frac{2n}{3m})t^4+(\frac{10nx_1^2+10nx_2^2-2m}{3m})t^2+(\frac{5nx_1^4+5nx_2^4-3my_2^2}{3m}).
\end{equation}
 \\

\noindent We see that if
 \\

 \begin{equation}\label{15}
 (\frac{5nx_1^4+5nx_2^4-3my_2^2}{3m}),
 \end{equation}
\\

\noindent to be a square, say $q^2$, then we may transform the above quartic to a cubic elliptic curve. For every  arbitrary values of $n$, $m$, this is done by choosing appropriate values of $x_1$, $x_2$, $y_2$.
 \\

\noindent  As an example, if $m=85$, $n=6$, we may set $x_1=1$, $x_2=2$, $y_2=1$, ($q=1$).
\\

\noindent Then \eqref{15} becomes
\\

\begin{equation}\label{16}
v^2=\frac{4}{85}t^4+\frac{26}{51}t^2+1.
\end{equation}
 \\

\noindent With the inverse transformation
\begin{equation}\label{117}
 t=\frac{2(X+\frac{26}{51})}{Y},
 \end{equation}

and

\begin{equation}\label{18}
v=-1+\frac{t^2X}{2},
\end{equation}
\\

\noindent the corresponding cubic  elliptic curve is

 \begin{equation}\label{19}
 Y^2=X^3+\frac{26}{51}X^2-\frac{16}{85}X-\frac{416}{4335}.
 \end{equation}
\\

\noindent The rank of this elliptic curve is $2$ and its generators are the points

\noindent $P_1=(X',Y')=(\frac{2}{3},\frac{28}{51})$, and $P_2=(X'',Y'')=(\frac{20777}{21675},\frac{-1908281}{1842375})$.\\

\noindent  Because of this, the  above elliptic curve has infinitely many rational points and we may obtain infinitely many solutions for the Diophantine equation too. \\
 By using the point $P_1$, we get $(t,v)=(\frac{30}{7},\frac{251}{49})$, by calculating $X_i$, $Y_i$, from  the above relations and after  some simplifications  and canceling  the  denominators of  $X_i$, $Y_i$, we obtain a solution for the Diophantine equation as
 \\

$6.(1813^5+1127^5+2156^5+784^5)=85.(158123^3+(-14063)^3+88837^3+55223^3)$.
\\

\noindent By choosing the other points on the elliptic curve such as

\noindent $nP_1$, $nP_2$ ($n=2, 3$, $\cdots$ ) we obtain infinitely many solutions for the  Diophantine equation.
  \\

\noindent By taking $m=17$, $n=3$, in the \eqref{14}, we may choose appropriate values  $x_1=1$, $x_2=2$, $y_2=1$, ($q=2$).
\\

\noindent Then \eqref{14} becomes
\\

\begin{equation}\label{501}
v^2=\frac{2}{17}t^4+\frac{116}{51}t^2+4.
 \end{equation}
 \\

\noindent With the inverse transformation
\begin{equation}\label{21}
 t=\frac{4(X+\frac{116}{51})}{Y},
 \end{equation}

and

\begin{equation}\label{22}
v=-2+\frac{t^2X}{4},
\end{equation}

\noindent the corresponding cubic  elliptic curve is

 \begin{equation}\label{23}
 Y^2=X^3+\frac{116}{51}X^2-\frac{32}{17}X-\frac{3712}{867}.
 \end{equation}
\\

\noindent The rank of this elliptic curve is $1$ and its generator is the point

\noindent $P=(X,Y)=(4,\frac{160}{17})$. Because of this, the  above elliptic curve has infinitely many rational points and we may obtain infinitely many solutions for the Diophantine equation too. \\
 Since $P=(X,Y)=(4,\frac{160}{17})$, we get $(t,v)=(\frac{8}{3},\frac{46}{9})$, by calculating $X_i$, $Y_i$, from the above relations and after  some simplifications  and canceling  the  denominators of  $X_i$, $Y_i$, we obtain a solution for the Diophantine equation as
 \\

$3.(99^5+45^5+126^5+18^5)=17.(1890^3+(-594)^3+891^3+405^3)$.
\\

\noindent By choosing the other points on the elliptic curve such as

\noindent $nP$ ($n=2, 3$, $\cdots$ ) we obtain infinitely many solutions for the  Diophantine equation.
\\

\noindent The Sage software has been used for calculating  the rank of  the elliptic curves. (see \cite{7})
\\


\end{document}